# From unicellular fatgraphs to trees

Thomas J. X. Li, Christian M. Reidys[*]

*Biocomplexity Institute of Virginia Tech*
*1015 Life Sciences Circle, Blacksburg, VA 24061, USA*

*Department of Mathematics of Virginia Tech*
*225 Stanger Street, Blacksburg, VA 24061, USA*

**Abstract**

In this paper we study the minimum number of reversals needed to transform a unicellular fatgraph into a tree. We consider reversals acting on boundary components, having the natural interpretation as gluing, slicing or half-flipping of vertices. Our main result is an expression for the minimum number of reversals needed to transform a unicellular fatgraph to a plane tree. The expression involves the Euler genus of the fatgraph and an additional parameter, which counts the number of certain orientable blocks in the decomposition of the fatgraph. In the process we derive a constructive proof of how to decompose non-orientable, irreducible, unicellular fatgraphs into smaller fatgraphs of the same type or trivial fatgraphs, consisting of a single ribbon. We furthermore provide a detailed analysis how reversals affect the component-structure of the underlying fatgraphs. Our results generalize the Hannenhalli-Pevzner formula for the reversal distance of signed permutations.

*Keywords:* Fatgraph, Genus, Reversal, Tree, Irreducible

## 1. Introduction

A *fatgraph* is an embedding of a connected finite graph into a 2-dimensional, compact surface, such that each connected component of its complement is homeomorphic to a disc. Fatgraphs are also referred to as *maps* or *ribbon graphs*, and are a special type of *dessin d'enfant* as defined by Grothendieck [9]. In this paper we consider fatgraphs of orientable as well as non-orientable surfaces. A fatgraph is *unicellular* if it has only one boundary component.

---
[*]Corresponding author.
*Email addresses:* `thomasli@vt.edu` (Thomas J. X. Li), `duck@santafe.edu` (Christian M. Reidys)

A well-studied example are plane trees, i.e. unicellular fatgraphs embedded into the sphere. Unicellular fatgraphs appear frequently in combinatorics [16, 15, 5, 6] and play a central role in geometry: moduli spaces of curves [11, 18], in knot theory: Bollobás-Riordan-Tutte polynomial [4, 7], in mathematical physics: quantum field theory [2, 14], and in bioinformatics: modeling of macromolecules [19] including proteins [20] as well as RNA structures [17, 1].

Recently, Huang and Reidys [13] give a topological interpretation of the Hannenhalli-Pevzner theory for the reversal distance of signed permutations. The reversal distance problem originates from genome rearrangements in molecular evolution and was analyzed in a seminal paper by Hannenhalli and Pevzner [10] employing breakpoint graphs. The topological framework of [13] constructs a bijection between signed permutations and $\pi$-maps, an equivalence class of particular unicellular fatgraphs, see Section 2 for details. Reversals on signed permutations correspond to three operations on $\pi$-maps, namely gluing, slicing and half-flipping. Fig. 1 illustrates *reversals* acting on an arbitrary unicellular fatgraph. In the following we will consider gluing, slicing and half-flipping on arbitrary fatgraphs and refer to these simply as reversals. The reversal distance of a signed permutation is connected to the topological genus of the associated $\pi$-map and can be formulated as the solution of a combinatorial optimization problem, namely as the minimum number of reversal actions needed to reduce a $\pi$-map to a tree.

In this paper we study the *r-distance* of a unicellular fatgraph, that is, the minimum number of reversals needed to transform it to a plane tree.

The main result of this paper is an explicit formula for computing the $r$-distance for arbitrary fatgraphs. The formula generalizes the Hannenhalli-Pevzner expression for the reversal distance of signed permutations and can be computed in polynomial time:

$$d(\mathbb{G}) = \begin{cases} g + h + 1 & \text{if } h \neq 1, h \text{ is odd and all } h \text{ E-blocks are super-blocks,} \\ g + h & \text{otherwise,} \end{cases}$$
(1.1)

where $g$ is the Euler genus of $\mathbb{G}$ and $h$ is the number of E-blocks in $\mathbb{G}$, see Section 9.

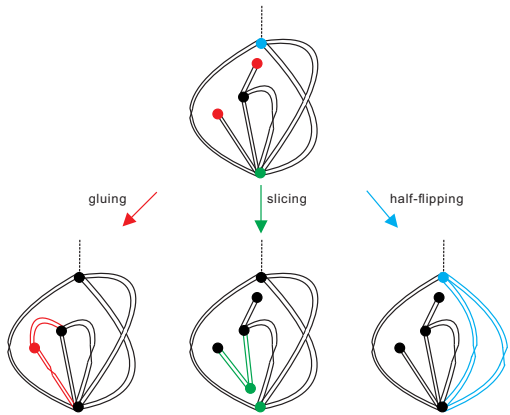

**Fig.** 1: The three actions on fatgraphs: gluing, slicing and half-flipping.

To prove eq. (1.1), we decompose a unicellular fatgraph into *components*, i.e. equivalence classes of associated ribbons. Given a unicellular fatgraph, two ribbons $r_1, r_2$ are *associated* if there exists a sequence of ribbons $(r_1 = w_1, w_2, \ldots, w_{k-1}, w_k = r_2)$, such that each two consecutive ribbons are traversed alternately on the boundary component, see Section 3. A component is *trivial* if it consists of a single bi-directional ribbon[1]. This decomposition allows us to introduce the *component tree*, $T$, see Section 3.

By deleting trivial components, $\mathbb{G}$ disects into *blocks*. This gives rise to the *block tree* of $\mathbb{G}$, see Section 3. A block containing only orientable components is called *orientable* and *exposed* (*E-block*) if it is orientable and not contained in a path in the block tree, joining two other orientable blocks.

Eq. (1.1) follows in three steps: first we show that any non-orientable component can be reduced into a tree via $g$ reversals. Here we can identify upfront the set of mono-directional ribbons, whose slicing produces a fatgraph consisting of exclusively non-orientable as well as trivial components. Compared with [13] and [10], our proof provides new insight and is constructive. We can explicitly characterize the relevant ribbons as opposed to showing the mere their existence via a proof by contradiction.

Calling a block, containing at least one non-orientable component, *non-orientable* we secondly show that any non-orientable block can be transformed into a tree via exactly $g$ slicings. The idea will be to successively merge two non-orientable and orientable components into a single, non-orientable component via slicings. By successively eliminating orientable components, we thereby reduce the problem to non-orientable components.

---

[1]corresponding to the cut edges, considered in graph theory

Thirdly, it remains to consider blocks containing exclusively orientable components. Their analysis involves certain paths within the block tree and blocks, not contained in these paths, i.e. *exposed* blocks (E-blocks). The E-blocks are of central importance for deriving eq. (1.1). Certain blocks contained in the above mentioned paths are merged into a single non-orientable block by means of gluings and half-flippings. This reduces the problem to non-orientable blocks.

Fatgraphs appear in the context of RNA pseudoknots, RNA structures that exhibit cross-serial interactions. Pseudoknot RNA is of relevance in the context of aging (RNA telomerase) [3] and disease (non-coding RNA) [8]. The topological classification of RNA pseudoknots, employing orientable fatgraphs, has been studied in [17, 1] and led to folding algorithms that construct a partition function of such configurations in $O(n^6)$ time and $O(n^4)$ space complexity [21]. [12] presents a context-free grammar for uniformly generating RNA structures of fixed topological genus. The grammar is facilitated by a bijection reducing the genus of orientable, unicellular fatgraphs by successively tri-slicing vertices [5].

The framework presented here deals with arbitrary fatgraphs and employs reversals, i.e. bi-slicing vertices, by which we transform arbitrary, unicellular fatgraphs into trees. We present a simple formula for the minimum number of reversals needed to obtain such a tree. This constitutes as a first step towards a bijection between unicellular fatgraphs and unicellular fatgraphs of lower genera along the lines of [5]. Such a "constructive" genus reduction has the potential to allow to recursively construct more complex RNA structures or proteins, extending the results of [12].

The paper is organized as follows: in Section 2, we introduce fatgraphs. In Section 3, we discuss irreducibility and components and introduce the notion of crossing ribbons and show that any unicellular fatgraph can be decomposed into components. In Section 4 we study reversals on fatgraphs and in Section 5 we analyze the effect of reversals on the direction and crossing of ribbons. We prove in Section 6 that for any non-orientable component its $r$-distance equals its genus $g$. Section 7 investigates the effect of reversals on the component tree $T$ and the block tree $B$. In Section 8, we show that any fatgraph with only non-orientable blocks has its $r$-distance $g$. Finally we analyze orientable blocks and prove an analogue of the Hannenhalli-Pevzner formula for the $r$-distance of arbitrary fatgraph in Section 9.

## 2. Fatgraphs

Following [13], a fatgraph can be expressed using sectors, i.e. oriented wedges, represented as pairs $(i, \omega(i))$, where $i$ is a label and $\omega(i) \in \{-1, +1\}$ is the orientation. We denote counterclockwise and clockwise orientations by $\omega(i) = +1$ and $\omega(i) = -1$, respectively. By abuse of notation, we shall identify a sector with its label and only consider a sector as a pair, when its orientation is of particular relevance.

**Definition 1.** A *rooted fatgraph* with $n$ ribbons is a quadruple $\mathbb{G} = ([2n + 1], \sigma, \gamma, \omega)$, where $[2n + 1] = \{1, \ldots, 2n + 1\}$ is a set of labeled sectors, $\sigma$ and $\gamma$ are permutations representing vertices and boundary components with $\gamma(2n+1) = \sigma(2n+1) = 1$, the function $\omega : [2n+1] \longrightarrow \{-1, +1\}$ represents the orientation of sectors with $\omega(1) = \omega(2n + 1)$, and finally for each pair $(x, \sigma(x))$ (where $x \neq 2n + 1$), there exists a unique pair $(y, \sigma(y))$ such that

- Case 1:
$$\begin{array}{ccc} x & \xrightarrow{\sigma} & \sigma(x) \\ \gamma \downarrow & & \downarrow \gamma \\ \sigma(y) & \xleftarrow{\sigma} & y \end{array}$$
with $\omega(x) = \omega(\sigma(y)), \omega(\sigma(x)) = \omega(y)$,

- Case 2:
$$\begin{array}{ccc} x & \xrightarrow{\sigma} & \sigma(x) \\ & \gamma \times & \\ \sigma(y) & \xleftarrow{\sigma} & y \end{array}$$
with $\omega(x) = -\omega(y), \omega(\sigma(x)) = -\omega(\sigma(y))$.

The directions of the $\gamma$-verticals are implied by the orientations of pairs of sectors $(x, \sigma(x))$ and $(y, \sigma(y))$ and we shall refer to the above diagrams as *untwisted* and *twisted* ribbons, respectively.

The above definition is equivalent to that in [13], since a rooted fatgraph can be obtained from a fatgraph via bisecting the sector, 1, into two sectors 1 and $2n + 1$, see Fig. 2 (A) and (B). In the presentation of a fatgraph, we will use an arrow or a sign to denote orientation.

A ribbon is given by two pairs of sectors, $((x, \sigma(x)), (y, \sigma(y)))$. Ribbons with mono- and bi-directional verticals are called $m$- *and* $b$-*ribbons*, respectively. For $m$-ribbons we have $\omega(x) = -\omega(\sigma(x))$ and $\omega(y) = -\omega(\sigma(y))$, while for $b$-ribbons $\omega(x) = \omega(\sigma(x))$ and $\omega(y) = \omega(\sigma(y))$ hold. An untwisted ribbon (also for twisted ribbon) can be either a $m$- or $b$-ribbon, see Fig. 2 (B).

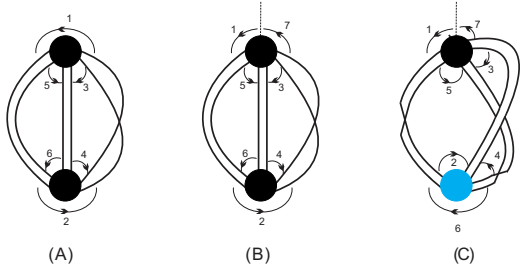

**Fig. 2**: Fatgraphs. (A) A non-orientable fatgraph as in [13]. (B) The same fatgraph, using our notation: $\sigma = (1,5,3,7)(2,4,6)$, $\gamma = (1,2,3,4,7)(5,6)$ and orientations $\omega(1) = \omega(2) = \omega(5) = \omega(6) = \omega(7) = +1$, $\omega(3) = \omega(4) = -1$. The ribbon $((1,5),(6,2))$ is untwisted and bi-directional, while the ribbon $((5,3),(4,6))$ is untwisted and mono-directional. (C) Flipping: vertex $(2,4,6)$ in (B) is flipped into $(2,6,4)$ in (C).

A fatgraph $\mathbb{G}$ represents a cell-complex of a surface $F(\mathbb{G})$: the topological quotient space $F(\mathbb{G})$ is obtained by identifying the sides of ribbons consistent with their orientations, see [20]. We call a fatgraph *orientable* if it induces an orientable surface $F(\mathbb{G})$, and *non-orientable*, otherwise. Since we work with orientable as well as non-orientable fatgraphs in this paper, we define the *genus* $g$ of $\mathbb{G}$ as the *Euler genus* $g^*$ (or *modified genus*) of $F(\mathbb{G})$, i.e.,

$$2 - g^* = v - e + b,$$

where $v, e, b$ are the numbers of vertices, ribbons and boundary components.

A *flip* of a vertex $v$ reverses the cyclic ordering of its sectors and changes their respective orientations, see Fig. 2 (B) and (C). Two fatgraphs $\mathbb{G}_1$ and $\mathbb{G}_2$ are *isomorphic* if $\mathbb{G}_2$ can be obtained from $\mathbb{G}_1$ by relabeling of sectors and flipping of vertices. By abuse of notion we shall refer to an equivalence class of isomorphic fatgraphs simply as a fatgraph.

A *unicellular* fatgraph is a fatgraph in which $\gamma$ is a unique cycle. Without loss of generality we can assume that $\gamma = (1, \ldots, 2n+1)$. In particular, a unicellular fatgraph of genus zero is a plane tree. A *$\pi$-map* [13] is a unicellular fatgraph having a distinguished vertex (called the center) and any other vertices (called external vertices), such that any ribbon is incident to the center and an external vertex. We generalize a criterion for orientability of $\pi$-maps [13] to arbitrary unicellular fatgraphs.

**Proposition 1** (Orientability). *A unicellular fatgraph $\mathbb{G}$ is non-orientable if and only if $\mathbb{G}$ contains a m-ribbon.*

## 3. Some basic facts

In this section, we introduce irreducibility and decompose unicellular fatgraphs into their components, extending the results of [13].

Let $\mathbb{G} = ([2n+1], \sigma, \gamma, \omega)$ be a unicellular fatgraph. Given a ribbon $r = ((x, \sigma(x)), (y, \sigma(y)))$, the *origin* $r^L$ and the *terminus* $r^R$ of $r$ are $r^L = \min_\gamma\{x, \sigma(x), y, \sigma(y)\}$ and $r^R = \max_\gamma\{x, \sigma(x), y, \sigma(y)\}$, respectively. In the traversal of the boundary component, $r^L$ is the first $r$-sector, whereas $r^R$ is the last.

Two ribbons $r_1$ and $r_2$ are $\mathbb{G}$-*crossing* or simply *crossing* if

$$r_1^L <_\gamma r_2^L <_\gamma r_1^R <_\gamma r_2^R \quad \text{or} \quad r_2^L <_\gamma r_1^L <_\gamma r_2^R <_\gamma r_1^R.$$

Otherwise, $r_1$ and $r_2$ are called $\mathbb{G}$-*non-crossing*. Two ribbons are $\mathbb{G}$-*associated* or *associated* if there exists a sequence of ribbons $(r_1 = w_1, w_2, \ldots, w_{k-1}, w_k = r_2)$ such that $w_i, w_{i+1}$ are $\mathbb{G}$-crossing.

**Definition 2. (Irreducibility)** A unicellular fatgraph $\mathbb{G}$ is *irreducible* if any two $\mathbb{G}$-ribbons $r_1$ and $r_2$ are associated.

Clearly, ribbon association is an equivalence relation and the set of ribbons partitions into equivalence classes of associated ribbons, see Fig. 3 (A).

**Definition 3. (Component)** A $\mathbb{G}$-*component* is an equivalence class of associated ribbons in $\mathbb{G}$. A $\mathbb{G}$-component is *trivial* if it consists of only one $b$-ribbon.

The following proposition is straightforward, see the Supplementary Material (SM).

**Proposition 2.** *For any unicellular fatgraph $\mathbb{G}$ the following assertions hold:*
*(1) $\mathbb{G}$ uniquely decomposes into a set of components $C$,*
*(2) each component $C$ induces an irreducible unicellular fatgraph $\mathbb{G}_C$,*
*(3) $g_\mathbb{G} = \sum_C g_{\mathbb{G}^C}$.*

Fig. 3 (A) depicts a decomposition of a fatgraph into components.

For any two sectors $i$ and $j$, we call $[i,j]_\gamma = \{s \colon i \leq_\gamma s \leq_\gamma j\}$ the *interval of $\gamma$* and $(i,j)_\gamma$ and $(i,j]_\gamma$ are defined analogously. The $C$-*trace* $I_C$ is the set of sectors of $C$-ribbons. Clearly, we have $I_C = \dot\cup_i^l [a_i, c_i]_\gamma$, where $a_1 <_\gamma \cdots <_\gamma a_l$ and $[c_i, a_{i+1}]_\gamma$ is a $C$-*gap*, where $1 \leq i \leq l$.

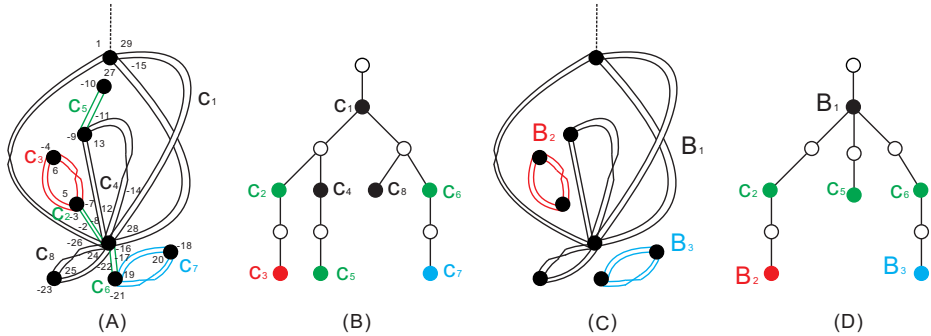

**Fig. 3**: (A) The decomposition of a unicellular fatgraph into its components. A unicellular fatgraph $\mathbb{G}$ of genus 6 is decomposed into the eight components $C_1, \ldots, C_8$. The latter have the genera $2, 0, 1, 1, 0, 0, 1, 1$, respectively. (B) The component tree $T$ of $\mathbb{G}$. (C) Deleting trivial components splits $\mathbb{G}$ into its blocks. (D) The block tree $B$ of $\mathbb{G}$.

**Proposition 3.** *Any two traces of components are either subsequent or nested, i.e. one is contained in a gap of the other.*

The proof is straightforward and given in the SM. Let $\mathbb{G}$ be a fatgraph having the components $C_1, \ldots, C_k$ with traces $I_{C_h} = \cup_i^{l_h} [a_i^{(h)}, c_i^{(h)}]_\gamma$, where $[c_i^{(h)}, a_{i+1}^{(h)}]_\gamma$ are $C_h$-gaps, $1 \leq h \leq k$. Clearly, the set of gaps is partially ordered by set inclusion.

**Definition 4. (Component Tree)** The *component tree* $T$ of $\mathbb{G}$ is a rooted, bicolored tree having $V_1 = \{C_1, \ldots, C_k\}$ (black) and $V_2 = \{[c_i^{(h)}, a_{i+1}^{(h)}]_\gamma \mid 1 \leq h \leq k, i\} \cup \{[1, 2n+1]_\gamma\}$ (white) as vertices. $T$ has the following adjacency relations: a white vertex is a child of a black vertex $C_h$ if it is a $C_h$-gap, a black vertex is a child of the white root if its trace is not contained in any gap, and a black vertex, $C_j$ is a child of a white vertex, $[c_i^{(h)}, a_{i+1}^{(h)}]_\gamma$, if $[c_i^{(h)}, a_{i+1}^{(h)}]_\gamma$ is the minimum gap containing the $C_j$-trace.

**Remark:** Proposition 3 guarantees that $T$ has no cycles and thus is well-defined.

Fig. 3 (A) and (B) present a unicellular fatgraph and its component tree.

Each trivial component corresponds to a cut vertex or a leave in the component tree $T$. By deleting all vertices corresponding to trivial components, we obtain a set of subtrees $B_h$ of $T$. A *block* $\mathbb{B}_h$ of $\mathbb{G}$ is the fatgraph induced by the components corresponding to $B_h$, see Fig. 3 (C).

The set of sectors of $\mathbb{B}_h$-ribbons is called the $\mathbb{B}_h$-*trace*, $I_{\mathbb{B}_h}$. The trace $I_{\mathbb{B}_h} = \cup_j^{k_h}[a_j^{(h)}, c_j^{(h)}]_\gamma$, where $[c_j^{(h)}, a_{j+1}^{(h)}]_\gamma$ are the $\mathbb{B}_h$-*gaps* for $1 \leq j < k_h$. As in Proposition 3, any two blocks are either subsequent around the boundary component or one is contained in a gap of the other. This property allows us to define a tree on blocks.

**Definition 5. (Block Tree)** The *block tree*, $B$, is the rooted bicolored tree having the set of trivial components and blocks, $V_1$ (black) and the set of all gaps of trivial components and blocks, as well as the root $[1, 2n+1]_\gamma$, $V_2$ (white) as vertex set. The adjacency relations in $T$ are as follows: a white vertex $v_2$ is a child of a black vertex $v_1$ if $v_2$ is a $v_1$-gap. A black vertex is a child of the white root if its trace is not contained in any gap. A black vertex $v_1$ is a child of a white vertex $v_2$ if $v_2$ is the minimal gap containing the $v_1$-trace.

**Remark:** A gap is called a *boundary gap* if it is contained in some $B_h$ and adjacent to some trivial component. Given a component tree $T$, by definition, its corresponding block tree $B$ is constructed by first contracting each subtree $B_h$ excluding boundary gaps into a black vertex $v_i$ and then each boundary gap is connected with its corresponding black vertex $v_i$ and trivial components, see Fig. 3 (B) and (D).

## 4. Reversals

In this section, we generalize the reversal operations on $\pi$-maps [13] to unicellular fatgraphs. Given two sectors $i, j$ with $i <_\gamma j$, we consider

**Gluing:** suppose $i$ and $j$ are located at two distinct vertices $v_1 = (i, i_1, \ldots, i_p)$ and $v_2 = (j, j_1 \ldots, j_q)$. Without loss of generality (flipping), we can always assume that $i$ and $j$ have different orientations, i.e. $\omega(i) = -\omega(j)$. Then gluing $v_1$ and $v_2$ along sectors $i$ and $j$, produces the new vertex

$$v = (i, j_1 \ldots, j_q, j, i_1, \ldots, i_p),$$

and the fatgraph $\tilde{\mathbb{G}}$, in which $v_1$ and $v_2$ merge into $v$, see Fig. 4.

**Slicing:** suppose $i$ and $j$ are located at $v = (i, j_1 \ldots, j_q, j, i_1, \ldots, i_p)$ and have different orientations, i.e. $\omega(i) = -\omega(j)$. Then by slicing $v$ along sectors $i$ and $j$, we derive the vertices $v_1$ and $v_2$

$$v_1 = (i, i_1, \ldots, i_p), \qquad v_2 = (j, j_1 \ldots, j_q),$$

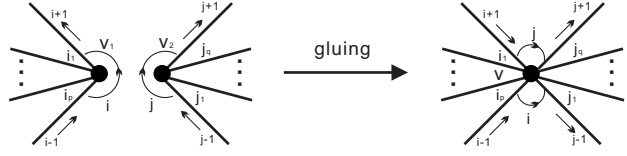

**Fig.** 4: Gluing.

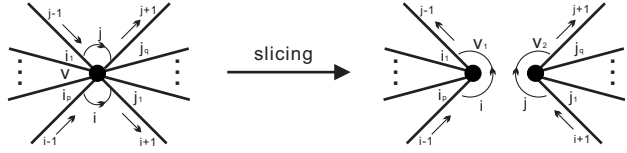

**Fig.** 5: Slicing.

and the fatgraph $\tilde{\mathbb{G}}$, in which $v$ is replaced by $v_1$ and $v_2$, see Fig. 5.

**Half-flipping:** suppose $i$ and $j$ are located at $v = (i, j_1 \ldots, j_q, j, i_1, \ldots, i_p)$, have the same orientation and $\omega(i) = \omega(j) = +1$. Then half-flipping is obtained as follows: we first split the vertex $v$ along the sectors $i$ and $j$ into $v_1$ and $v_2$, secondly flip $v_2$ containing the sector $j$ and thirdly glue $v_1$ and $v_2$ along sectors $i$ and $j$ into $\tilde{v}$, see Fig. 6. This sequence of operations is tantamount to flipping all ribbons attaching to $v$ from sector $i$ to sector $j$ counterclockwise (clockwise if $\omega(i) = \omega(j) = -1$) and produces the new fatgraph $\tilde{\mathbb{G}}$, in which $v$ is replaced by $\tilde{v} = (i, j_q \ldots, j_1, j, i_1, \ldots, i_p)$, see Fig. 6.

**Lemma 1.** *Let $\mathbb{G}$ be a unicellular fatgraph with the boundary component*

$$\gamma = (\underline{1, \ldots, i-1, i, i+1, \ldots, j-1, j, j+1, \ldots, 2n+1}).$$

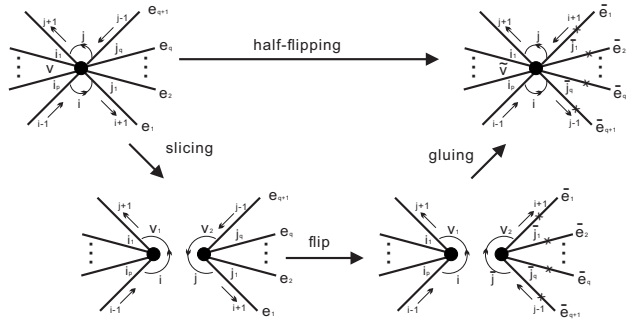

**Fig.** 6: Half-flipping.

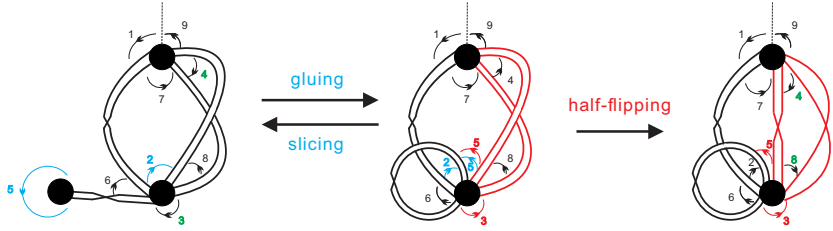

**Fig. 7**: Three types of reversals: gluing, slicing and half-flipping.

Then $\tilde{\mathbb{G}}$ is unicellular and has the boundary component

$$\tilde{\gamma} = (\underline{1,\ldots,i-1,i},\underline{j-1,\ldots,i+1},\underline{j,j+1,\ldots,2n+1}). \tag{4.1}$$

In case of gluing or slicing, we have

$$\tilde{\omega}(k) = \begin{cases} \omega(k) & \text{if } k \leq_\gamma i \text{ or } k \geq_\gamma j \\ -\omega(k) & \text{if } i+1 \leq_\gamma k \leq_\gamma j-1. \end{cases} \tag{4.2}$$

Thus gluing, slicing, or half-flipping manifest as a reversal of the boundary component and we shall refer to them as *reversals* or *i,j-reversals*. In Fig. 7, we illustrate the three types of reversals.

The *r-distance*, $d(\mathbb{G})$, of a unicellular fatgraph $\mathbb{G}$ is the minimum number of reversals needed to transform it into a unicellular fatgraph of genus zero, i.e. a tree.

**Corollary 1.** *Any unicellular fatgraph $\mathbb{G}$ of genus $g$ has $d(\mathbb{G}) \geq g$.*

By construction, each reversal reduces the number of vertices by at most one. The corollary follows directly from Euler's characteristic formula.

### 5. Reversals and ribbons

In this section, we investigate how reversals affect ribbons. Let $r$ denote a $\mathbb{G}$-ribbon and $r^{i,j}$ denote the $\tilde{\mathbb{G}}$-ribbon that is produced from $r$ by the $i,j$-reversal. By construction, there is a natural bijection between ribbons of $\mathbb{G}$ and those of $\tilde{\mathbb{G}}$.

A ribbon $r$ *intersects* an interval $[i,j]_\gamma$ if $r^L <_\gamma i <_\gamma r^R \leq_\gamma j$ or $i \leq_\gamma r^L <_\gamma j <_\gamma r^R$. In some sense, this generalizes the notion of crossing: two ribbons $r_1$ and $r_2$ cross in $\mathbb{G}$ if and only if $r_1$ intersects the interval $[r_2^L, r_2^R]_\gamma$. The *directional status* of a ribbon $r$ indicates it is mono- or bi-directional.

**Lemma 2.** *The directional status of $r$ changes after the $i,j$-reversal if and only if $r$ intersects the interval $[i,j]_\gamma$.*

The lemma follows via direct verification of the orientation of the two pairs $\{r^L, \gamma(r^L)\}$ and $\{\gamma^{-1}(r^R), r^R\}$, utilizing eq. (4.2), see the SM.

Let $r_1, r_2$ be two ribbons and let $\text{cr}(r_1, r_2) = 1, 0$ depending on whether $r_1, r_2$ are crossing or not. We call $\text{cr}(r_1, r_2)$ the *crossing status*.

**Lemma 3.** *An $i,j$-reversal changes $cr(r_1, r_2)$ if and only if $r_1, r_2$ intersect $[i,j]_\gamma$.*

*Proof.* We first establish a criterion for maintaining the crossing status of two ribbons.

**Criterion:** if we have the equality of sets $(r_1^L, r_1^R)_\gamma = ((r_1^{i,j})^L, (r_1^{i,j})^R)_{\tilde{\gamma}}$, then $r_1$ and $r_2$ cross in $\mathbb{G}$ if and only if $r_1^{i,j}$ and $r_2^{i,j}$ cross in $\tilde{\mathbb{G}}$.

In view of $(r_1^L, r_1^R)_\gamma = ((r_1^{i,j})^L, (r_1^{i,j})^R)_{\tilde{\gamma}}$, the set of intermediate sectors between origin and terminus of $r_1$ and $r_1^{i,j}$ remains unchanged. We verify that $r_2$ intersects $[r_1^L, r_1^R]_\gamma$ if and only if $r_2^{i,j}$ intersects $[(r_1^{i,j})^L, (r_1^{i,j})^R]_{\tilde{\gamma}}$.

**Claim 1:** if $r_1$ does not intersect the interval $[i,j]_\gamma$, then $r_1$ and $r_2$ cross in $\mathbb{G}$ if and only if $r_1^{i,j}$ and $r_2^{i,j}$ cross in $\tilde{\mathbb{G}}$.

We distinguish the four cases: $r_1^L <_\gamma i <_\gamma j <_\gamma r_1^R$, $i \leq_\gamma r_1^L <_\gamma r_1^R \leq_\gamma j$, $r_1^L <_\gamma r_1^R \leq_\gamma i <_\gamma j$ and $i <_\gamma j \leq_\gamma r_1^L <_\gamma r_1^R$. Since all cases can be argued analogously, it suffices to prove the first: in view of eq. (4.1), the ribbon $r_1^{i,j}$ has the same origin and terminus as $r_1$, i.e. $(r_1^{i,j})^L = r_1^L$ and $(r_1^{i,j})^R = r_1^R$. Moreover, $(r_1^{i,j})^L <_{\tilde{\gamma}} i <_{\tilde{\gamma}} j <_{\tilde{\gamma}} (r_1^{i,j})^R$. The boundary component $\tilde{\gamma}$ reverses only the sequence $i+1, \ldots, j-1$, contained in $[r_1^L, r_1^R]_\gamma$ and $[(r_1^{i,j})^L, (r_1^{i,j})^R]_{\tilde{\gamma}}$. Thus $[r_1^L, r_1^R]_\gamma = [(r_1^{i,j})^L, (r_1^{i,j})^R]_{\tilde{\gamma}}$ and the Criterion implies that $r_1, r_2$ cross in $\mathbb{G}$ if and only if $r_1^{i,j}, r_2^{i,j}$ cross in $\tilde{\mathbb{G}}$.

It remains to show

**Claim 2:** if $r_1$ and $r_2$ intersect the interval $[i,j]_\gamma$, then they change their crossing status.

Without loss of generality we may assume that $r_1^L <_\gamma r_2^L$ and $r_1, r_2$ are $\mathbb{G}$-crossing. We distinguish three cases: $i \leq_\gamma r_1^L <_\gamma r_2^L <_\gamma j <_\gamma r_1^R <_\gamma r_2^R$, $r_1^L <_\gamma i \leq_\gamma r_2^L <_\gamma r_1^R \leq_\gamma j <_\gamma r_2^R$ and $r_1^L <_\gamma r_2^L <_\gamma i <_\gamma r_1^R <_\gamma r_2^R \leq_\gamma j$. It suffices to check the first case, as the other two follow analogously. Then we have

$$i \leq_\gamma r_1^L <_\gamma \gamma(r_1^L) \leq_\gamma r_2^L <_\gamma \gamma(r_2^L) \leq_\gamma j \leq_\gamma \gamma^{-1}(r_1^R) <_\gamma r_1^R \leq_\gamma \gamma^{-1}(r_2^R) <_\gamma r_2^R.$$

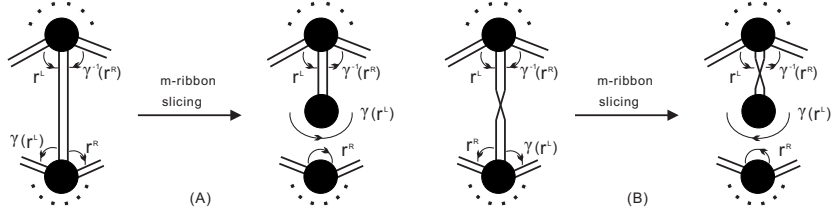

**Fig.** 8: (A) Slicing an untwisted $m$-ribbon. (B) Slicing a twisted $m$-ribbon.

Since the boundary component $\tilde{\gamma}$ differs only by reversing the part $i+1,\ldots,j-1$, we obtain

$$i \leq_{\tilde{\gamma}} \gamma(r_2^L) <_{\tilde{\gamma}} r_2^L \leq_{\tilde{\gamma}} \gamma(r_1^L) <_{\tilde{\gamma}} r_1^L \leq_{\tilde{\gamma}} j \leq_{\tilde{\gamma}} \gamma^{-1}(r_1^R) <_{\tilde{\gamma}} r_1^R \leq_{\tilde{\gamma}} \gamma^{-1}(r_2^R) <_{\tilde{\gamma}} r_2^R,$$

and $(r_1^{i,j})^L = \gamma(r_1^L)$, $(r_1^{i,j})^R = r_1^R$, $(r_2^{i,j})^L = \gamma(r_2^L)$, $(r_2^{i,j})^R = r_2^R$. Thus

$$i \leq_{\tilde{\gamma}} (r_2^{i,j})^L <_{\tilde{\gamma}} (r_1^{i,j})^L <_{\tilde{\gamma}} j <_{\tilde{\gamma}} (r_1^{i,j})^R <_{\tilde{\gamma}} (r_2^{i,j})^R,$$

i.e. $r_1^{i,j}$ and $r_2^{i,j}$ do not cross in $\tilde{\mathbb{G}}$. □

**Remark:** In the following we oftentimes write $r$ instead of $r^{i,j}$.

## 6. Non-orientable irreducible fatgraphs

In this section we show that any non-orientable component $C_i$ of $\mathbb{G}$ and thus any non-orientable, irreducible fatgraph, can be transformed into a tree using $g_i$ reversals.

We begin by considering a particular class of slicings, namely that of $m$-ribbons, $r$, at $\gamma(r^L)$ and $r^R$. Any $m$-ribbon $r$ consists of a set of four sectors, $\{r^L, \gamma(r^L), \gamma^{-1}(r^R), r^R\}$, where $\gamma(r^L)$ and $r^R$ ($r^L$ and $\gamma^{-1}(r^R)$) are located at the same vertex and have different orientations, see Fig. 8. Thus, slicing a $m$-ribbon is well-defined and results in a new, unicellular fatgraph $\tilde{\mathbb{G}}$. Slicing $r$ produces the set $\{r^L, \gamma(r^L), \gamma^{-1}(r^R)\}$, i.e. $r$ becomes a trivial component of $\tilde{\mathbb{G}}$, see Fig. 8.

Since $r_1$ and $r$ are $\mathbb{G}$-crossing if and only if $r_1$ intersects $[\gamma(r^L), r^R]_\gamma$, we derive the following corollary of Lemma 3.

**Corollary 2.** *Suppose that $r_1, r_2 \neq r$. Then slicing $r$, changes the crossing status of $r_1$ and $r_2$ if and only if $r_1$ and $r_2$ cross $r$.*

**Lemma 4.** *Any irreducible, non-orientable unicellular fatgraph, $\mathbb{G}$, can be sliced via a m-ribbon into non-orientable and trivial components.*

This is the analogue of Theorem 4 of [10] on signed permutations and Lemma 10 of [13] on $\pi$-maps. The proof given here is new: we identify upfront a distinguished $m$-ribbon, whose slicing facilitates the lemma.

*Proof.* For a fixed $m$-ribbon $r$, any other ribbon can be bi- (b) or mono-directional (m) and either cross $r$ (+) or not ($-$). This partitions the set of ribbons into $E_b^+(r), E_b^-(r), E_m^+(r)$ and $E_m^-(r)$. Then $z = |E_b^+(r)| + |E_m^-(r)|$ is a function of $r$. Let $e$ be a $m$-ribbon that maximizes $z$ and $\mathbb{G}^e$ be the fatgraph obtained by slicing $e$.

To prove the lemma, it suffices to show that any ribbon, $r_1$, is contained in a non-orientable $\mathbb{G}^e$-component or a trivial $\mathbb{G}^e$-component. To this end we inspect how slicing $e$ affects the directional- and crossing status of $r_1$.

- Case 1: $r_1 \in E_b^+(e) \cup E_m^-(e)$. Lemma 2 shows that $r_1$ is mono-directional in $\mathbb{G}^e$, whence $r_1$ is contained in a non-orientable component.

- Case 2: $r_1 \in E_m^+(e)$. We distinguish two sub-cases:

  - Case 2(a): there exists a ribbon $r_2$, such that $r_2 \in E_b^+(e) \cup E_m^-(e)$ and $r_2 \notin E_b^+(r_1) \cup E_m^-(r_1)$.
    If $r_2 \in E_b^+(e)$, then $r_2$ is a $b$-ribbon of $\mathbb{G}$ and $r_2 \in E_b^-(r_1)$. Since $r_2$ and $r_1$ cross $e$, Corollary 2 implies that their crossing status changes from $\mathbb{G}$-non-crossing to $\mathbb{G}^e$-crossing. Case 1 shows that $r_2$ is a $m$-ribbon in $\mathbb{G}^e$ and thus $r_1$ is contained in a non-orientable component.
    If $r_2 \in E_m^-(e)$, then $r_2$ is a $m$-ribbon of $\mathbb{G}$ and $r_2 \in E_m^+(r_1)$. Since $r_2$ does not cross $e$, Corollary 2 guarantees that $r_2$ does not change its crossing status after slicing $e$. Since $r_2$ and $r_1$ are $\mathbb{G}$-crossing, they are $\mathbb{G}^e$-crossing, which implies that $r_1$ crosses $m$-ribbon $r_2$ in $\mathbb{G}^e$ and thus is in a non-orientable component.

  - Case 2(b): for $r_2 \in E_b^+(e) \cup E_m^-(e)$, we have $r_2 \in E_b^+(r_1) \cup E_m^-(r_1)$. This implies $E_b^+(e) \subset E_b^+(r_1)$ and $E_m^-(e) \subset E_m^-(r_1)$. Since $|E_b^+(e)| + |E_m^-(e)|$ is maximal, we have $E_b^+(e) = E_b^+(r_1)$ and $E_m^-(e) = E_m^-(r_1)$. Moreover $E_m^-(e) = E_b^-(r_1)$ and $E_m^+(e) \setminus \{r_1\} = E_m^+(r_1) \setminus \{e\}$. This means that $e$ and $r'$ are $\mathbb{G}$-crossing if and only if $r_1$ and $r'$ are $\mathbb{G}$-crossing for any other ribbon $r'$. By Corollary 2, $e$ and $r'$ are $\mathbb{G}^e$-crossing if and only if $r_1$ and $r'$ are $\mathbb{G}^e$-crossing. Therefore both $e$ and $r_1$ are trivial components of $\mathbb{G}^e$.

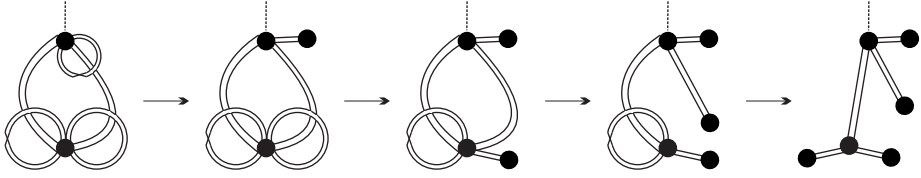

**Fig.** 9: From an irreducible, non-orientable unicellular fatgraph to a plane tree, by successively slicing four mono-directional ribbons.

- Case 3: $r_1 \in E_b^-(e)$. Since $\mathbb{G}$ is irreducible, there exists a sequence of ribbons $r_1, \ldots, r_k$ such that $r_i$ and $r_{i+1}$ are $\mathbb{G}$-crossing and $r_i \in E_b^-(e)$ for $1 \leq i \leq k-1$, but $r_k \notin E_b^-(e)$. We shall show that $r_k$ does not belong to Case 2(b): since $r_i$ does not cross $e$, Corollary 2 guarantees that $r_i, r_{i+1}$ are $\mathbb{G}^e$-crossing and $r_1, r_k$ are $\mathbb{G}^e$-associated, whence $r_k$ is not a trivial component. As a result $r_k$ belongs to Case 1 and Case 2(a) and is consequently associated with some $m$-ribbon in $\mathbb{G}^e$, and thus $r_1, r_k$ are in a non-orientable component.

□

**Theorem 1.** *Any irreducible, non-orientable fatgraph $\mathbb{G}$ has $d(\mathbb{G}) = g$.*

Fig. 9 depicts a sequence of slicings, that transforms an irreducible, non-orientable unicellular fatgraph into a tree.

## 7. Reversals on components and blocks

In this section we study $i, j$-reversals, mapping $\mathbb{G}$ into $\tilde{\mathbb{G}}$. Given the component tree $T$, we define a partial order on $T$-vertices by setting $v_i <_T v_j$ if $v_i$ is a descendant of $v_j$ in $T$. Suppose $[p, q]_\gamma$ is a $T$-gap and the set of $[p,q]_\gamma$-children is $\{C_h\}$ such that $I_{C_h} = [i^{(h-1)}, c_1^{(h)}]_\gamma \cup \cdots \cup [a_{l_h}^{(h)}, i^{(h)}]_\gamma$, where $p = i^{(0)} <_\gamma i^{(1)} <_\gamma \cdots <_\gamma i^{(k)} = q$.

**Definition 6.** A sector $s$ is *attached to a component $C$* if $s \in \cup_i^l (a_i, c_i)_\gamma$, where $I_C = \cup_i^l [a_i, c_i]_\gamma$. A sector $s$ is *attached to* $[p, q]_\gamma$ if $s \in \{i^{(0)}, i^{(1)}, \ldots, i^{(k-1)}, i^{(k)}\}$.

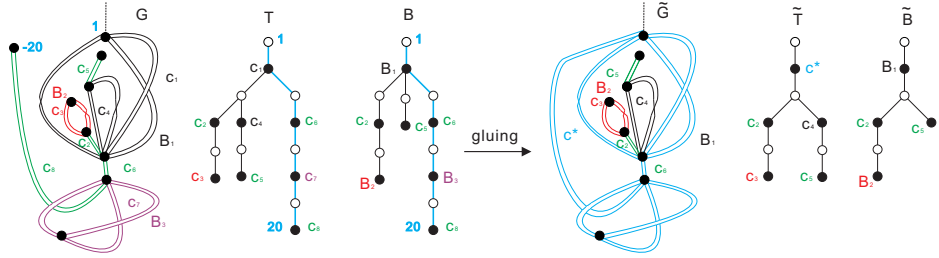

**Fig.** 10: The effect of gluing on the component tree and the block tree.

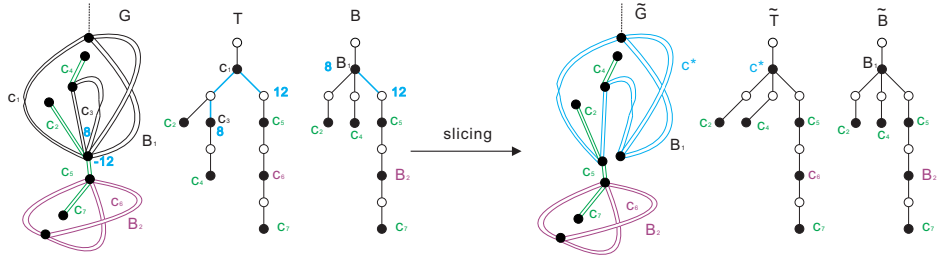

**Fig.** 11: The effect of slicing on the component tree and the block tree.

**Remark:** by construction, each sector is attached to a unique vertex in the component tree.

Since there is a natural bijection between ribbons of $\mathbb{G}$ and those of $\tilde{\mathbb{G}}$, we may, by abuse the notion, describe a $\tilde{\mathbb{G}}$-component or a $\tilde{\mathbb{G}}$-block by specifying its $\mathbb{G}$-ribbons. For example, we may say several $\mathbb{G}$-components form one $\tilde{\mathbb{G}}$-component, meaning that a $\tilde{\mathbb{G}}$-component consists of the ribbons of these $\mathbb{G}$-components.

Suppose $i, j$ are attached to $v_i$ and $v_j$ in the component tree, $T$. Let $P^{i,j}$ denote the $T$-path joining $v_i$ and $v_j$.

**Lemma 5. (Reversal Lemma)** *Each $\mathbb{G}$-component not on the path $P^{i,j}$ is a $\tilde{\mathbb{G}}$-component with its orientability unchanged. Furthermore, if the path $P^{i,j}$ contains at least two $\mathbb{G}$-components, then all components on $P^{i,j}$ form the $\tilde{\mathbb{G}}$-component $C_*$.*

We depict the effect of gluing, slicing and half-flipping for component trees in Fig. 10, Fig. 11 and Fig. 12, respectively.

*Proof.* We first prove two claims.

**Claim 1:** given a component $C$ with $I_C = \cup_i^l [a_i, c_i]_\gamma$, the following statements are equivalent

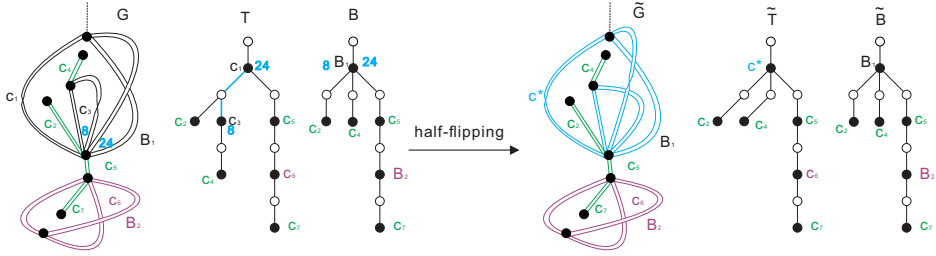

**Fig.** 12: The effect of half-flipping on the component tree and the block tree.

1. $a_1 <_\gamma i <_\gamma c_l$,
2. there exists a $C$-ribbon $r$ such that $r^L <_\gamma i <_\gamma r^R$,
3. $v_i \leq_T C$, i.e. vertex $v_i$ is either $C$ or a descendant of $C$ in $T$.

(1)$\Rightarrow$(2): suppose to the contrary that for each $C$-ribbon $r$, we have either $r^L \geq_\gamma i$ or $r^R \leq_\gamma i$. Note that any ribbon $r_1$ satisfying $r_1^L \geq_\gamma i$ does not cross with any ribbon $r_2$ with $r_2^R \leq_\gamma i$. Since any two $C$-ribbons are associated, we can assume that $r^L \geq_\gamma i$ holds for all $C$-ribbons. Thus $a_1 \geq_\gamma i$, a contradiction. (2)$\Rightarrow$(1) is straightforward. (1)$\Leftrightarrow$(3): by construction, if a gap $[s_1, s_2]_\gamma$ is a child of $C$, then $[s_1, s_2]_\gamma \subset [a_1, c_l]_\gamma$. Similarly, if $C$ is a child of a gap $[s_1, s_2]_\gamma$, then $[a_1, c_l]_\gamma \subset [s_1, s_2]_\gamma$. Therefore, we have that any $i \in (a_1, c_l)_\gamma$, if and only if $i$ is attached to $C$ or its descendant in $T$.

**Claim 2:** a component $C$ contains a ribbon $r$ intersecting the interval $[i, j]_\gamma$ if and only if $C$ is contained in $P^{i,j}$ ($P^{i,j}$-component) in $T$.

$\Leftarrow$: Suppose $C$ is a $P^{i,j}$-component, we have three cases

1. $v_i \leq_T C$ but $v_j \not\leq_T C$. By Claim 1, $a_1 <_\gamma i <_\gamma c_l \leq_\gamma j$ and there exists a $C$-ribbon $r$ such that $r^L <_\gamma i <_\gamma r^R \leq_\gamma j$, i.e. $r$ intersects the interval $[i, j]_\gamma$.
2. $v_j \leq_T C$ but $v_i \not\leq_T C$: this is analogous to the previous case.
3. $v_i \leq_T C$ and $v_j \leq_T C$. Suppose that any $C$-ribbon $r$ does not intersect the interval $[i, j]_\gamma$, that is, $r^R \leq_\gamma i$ or $j \leq_\gamma r^L$ or $r^L <_\gamma i <_\gamma j <_\gamma r^R$. Since the set of $C$-sectors consists of the sectors of its ribbons, $[i, j]_\gamma$ is contained in a $C$-gap. I.e. $i, j$ are contained in the same $C$-gap. As an $T$-ancestor of this gap, $C$ is not a $P^{i,j}$-component, a contradiction. Therefore there exists a $C$-ribbon $r$ intersecting the interval $[i, j]_\gamma$.

$\Rightarrow$: $r$ intersects the interval $[i, j]_\gamma$ and, without loss of generality, we can assume that $r^L <_\gamma i <_\gamma r^R \leq_\gamma j$. Claim 1 implies that $C$ is either $v_i$ or its ancestor in $T$. If $C$ is not a $P^{i,j}$-component, then $C$ is an ancestor of the minimum common ancestor of $v_i$ and $v_j$ in $T$. By construction, $i,j$ are contained in the same $C$-gap, whence $[i,j]_\gamma$ is not intersected by any $C$-ribbon, a contradiction to the assumption that a $C$-ribbon $r$ intersects the interval $[i,j]_\gamma$.

By Claim 2, each ribbon of a component not contained in $P^{i,j}$ does not intersect the interval $[i,j]_\gamma$. Lemma 3 guarantees that the reversal does not change its crossing status with any other ribbon. Therefore, each component not contained in $P^{i,j}$, is a $\tilde{\mathbb{G}}$-component with its orientability unchanged.

For each $P^{i,j}$-component, $C$, the set of ribbons $E_C$ splits into $E_C^+$ and $E_C^-$, representing the ribbons intersecting and not intersecting $[i,j]_\gamma$, respectively. By Claim 2, we know that $E_C^+ \neq \varnothing$. First suppose that $r_1, r_2$ are two ribbons intersecting the interval $[i,j]_\gamma$. Claim 2 guarantees that they belong to $P^{i,j}$-components. If $r_1, r_2$ belong to two different $P^{i,j}$-components, then $r_1, r_2$ are $\tilde{\mathbb{G}}$-crossing. This follows from Lemma 3 since both of them intersect the interval $[i,j]_\gamma$ and change their crossing status from $\mathbb{G}$-non-crossing to $\tilde{\mathbb{G}}$-crossing. If $r_1, r_2$ belong to the same $P^{i,j}$-component, $C$, we use that fact that $P^{i,j}$ contains at least two distinct components. In view of this fact, we can choose $r \in E_{C_1}^+$, where $C_1$ is a $P^{i,j}$-component other than $C$. Then $r_i$ and $r$ are $\tilde{\mathbb{G}}$-crossing for $i = 1, 2$, which implies that $r_1$ and $r_2$ are $\tilde{\mathbb{G}}$-associated. Secondly, for each ribbon $w_1$ in $E_C^-$, there exists a sequence of ribbons $(w_1, w_2, \ldots, w_{k-1}, w_k)$ such that $w_i, w_{i+1}$ are $\mathbb{G}$-crossing, $w_k \in E_C^+$ and $w_i \in E_C^-$ for $1 \leq i \leq k-1$. Lemma 3 guarantees that the reversal does not change the crossing status of a $E_C^-$-ribbon with any other ribbon, whence $w_i$ and $w_{i+1}$ are $\tilde{\mathbb{G}}$-crossing. Thus each ribbon $w_1 \in E_C^-$ is still associated to $w_k \in E_C^+$ in $\tilde{\mathbb{G}}$. Accordingly, any two ribbons of $P^{i,j}$-components become $\tilde{\mathbb{G}}$-associated, that is, all components on the path $P^{i,j}$ form one component. □

**Remark:** extending the proof of Lemma 5, we can construct the component tree $\tilde{T}$ of $\tilde{\mathbb{G}}$ utilizing Lemma 1, see the SM.

We call two components *adjacent* if their corresponding vertices in the component tree have distance two, i.e. either one is the child of a gap of the other component, or both are children of the same gap. Two components, that are adjacent and contain some $m$-ribbon incident to their common vertex, are called *m-adjacent*.

**Proposition 4.** *If $C$ is a $\mathbb{G}$-component not contained in $P^{i,j}$ that is adjacent to a $P^{i,j}$-gap or a $P^{i,j}$-component, then $C$ is adjacent to $C_*$ in $\tilde{\mathbb{G}}$.*

**Proposition 5.** *Any two adjacent $\mathbb{G}$-components not contained in $P^{i,j}$ are either again adjacent or both are adjacent to $C_*$ in $\tilde{\mathbb{G}}$. In particular, any two adjacent $\mathbb{G}$-components not contained in $P^{i,j}$ are in the same block of $\tilde{\mathbb{G}}$.*

We next analyze how reversals affect blocks. A block, $\mathbb{B}$, is *covered* by a path, $Q$, in the block tree if $\mathbb{B}$ is either contained in $Q$ or $\mathbb{B}$ is adjacent to a gap contained in $Q$.

Suppose $v_i$ and $v_j$ belong to blocks $B_i$ and $B_j$, respectively. Let $Q^{i,j}$ denote the path joining $\mathbb{B}_i$ and $\mathbb{B}_j$ in the block tree $B$ of $\mathbb{G}$.

**Lemma 6.** *Each $\mathbb{G}$-block not covered by $Q^{i,j}$ is a $\tilde{\mathbb{G}}$-block of the same orientability. Furthermore, if the path $P^{i,j}$ contains at least two $\mathbb{G}$-components, then all blocks covered by $Q^{i,j}$ together with the trivial components contained in $Q^{i,j}$ form the $\tilde{\mathbb{G}}$-block $B_*$.*

*Proof.* In the block tree, $B$, each $\mathbb{G}$-block not covered by $Q^{i,j}$ is adjacent to some trivial components not contained in $Q^{i,j}$. Proposition 5 implies that this block is adjacent to the same set of trivial components and thus consists of the same components in $\tilde{\mathbb{G}}$ and this $\tilde{\mathbb{G}}$-block has the same orientability. Suppose that $\mathbb{B}$ is a $\mathbb{G}$-block covered by $Q^{i,j}$. By definition, there exists a $\mathbb{B}$-component, $C_\mathbb{B}$, that is either contained in $Q^{i,j}$ or adjacent to a gap contained in $Q^{i,j}$. If $C_\mathbb{B}$ is contained in $Q^{i,j}$, Lemma 5 shows that $C_*$ contains all $C_\mathbb{B}$-ribbons. If $C_\mathbb{B}$ is adjacent to a $Q^{i,j}$-gap, Proposition 4 implies that $C_\mathbb{B}$ is adjacent to $C_*$ in $\tilde{\mathbb{G}}$. Proposition 5 guarantees that any $\mathbb{B}$-component is contained in the $\tilde{\mathbb{G}}$-block containing $C_*$, i.e. contained in $\mathbb{B}_*$. Therefore, $\mathbb{B}_*$ consists of ribbons of trivial $Q^{i,j}$-components and blocks covered by $Q^{i,j}$. □

## 8. The Pac-Man game

A block that contains at least one non-orientable component is called *non-orientable* and *orientable*, otherwise. A fatgraph whose blocks are all non-orientable is called *block-non-orientable*.

In this section, we shall show that block-non-orientable fatgraphs of genus $g$ have $r$-distance $g$.

**Proposition 6. (Pac-Man Proposition)** *Suppose we have an orientable $\mathbb{G}$-component, $C_1$, and a non-orientable $\mathbb{G}$-component, $C_2$, that are m-adjacent. Then there exists an $i_1, i_2$-slicing that generates a fatgraph $\tilde{\mathbb{G}}$ such that*
(1) *$C_1$ and $C_2$ merge into a non-orientable $\tilde{\mathbb{G}}$-component,*
(2) *any other $\mathbb{G}$-component is a $\tilde{\mathbb{G}}$-component of the same orientability.*
(3) *if $\mathbb{G}$ is block-non-orientable, then $\tilde{\mathbb{G}}$ is.*

*Proof.* We shall apply Lemma 5 and to that end identify two sectors $i_1, i_2$ such that $P^{i_1,i_2}$ contains $C_1$ and $C_2$ but no other components. Since $C_1$ and $C_2$ are adjacent, there exists a gap $[g_1, g_2]_\gamma$ adjacent to both $C_1$ and $C_2$ in $T$. Let $v$ be the common vertex of $C_1$ and $C_2$. Since neither $C_1$ nor $C_2$ are trivial components, each contains at least two ribbons incident to $v$. Consequently, there exist at least three $C_1$-sectors as well as three $C_2$-sectors at $v$. Let $i_1$ be a sector of $I_{C_1} \setminus \{g_1, g_2\}$ located at $v$. Since $C_1$ and $C_2$ are m-adjacent and all $C_1$-ribbons are b-ribbons, $C_2$ contains a m-ribbon incident to $v$. That is, there exist two $C_2$-sectors at $v$ with different orientations. Thus we can choose $i_2$ to be a sector of $I_{C_2} \setminus \{g_1, g_2\}$ at $v$ with an orientation, different from $i_1$. Since $i_k \in I_{C_k} \setminus \{g_1, g_2\}$, the sector $i_k$ is attached to either the component $C_k$ or a gap incident to $C_k$, where $k = 1, 2$. Furthermore this gap cannot coincide with $[g_1, g_2]_\gamma$. This implies that the only two $P^{i_1,i_2}$-components are $C_1$ and $C_2$.

In order to prove (1) and (2), Lemma 5 reduces the work to showing that $C_*$ is non-orientable. Claim 2 of Lemma 5 guarantees that $C_1$ contains a ribbon, $r$, intersecting the interval $[i, j]_\gamma$. By Lemma 2, $r$ changes from a b-ribbon in $\mathbb{G}$ to a m-ribbon of $C_*$ in $\tilde{\mathbb{G}}$, implying that $C_*$ is non-orientable. To establish (3), we let $\mathbb{B}_i$ be the block containing $C_1$ and $C_2$. By construction, the path $Q^{i,j}$ covers only the block $\mathbb{B}_i$. Lemma 6 implies, that each $\mathbb{G}$-block other than $\mathbb{B}_i$ remains a non-orientable $\tilde{\mathbb{G}}$-block and $\mathbb{B}_i$ becomes the $\tilde{\mathbb{G}}$-block $B_*$ containing $C_*$. Since $C_*$ is non-orientable, $B_*$ is non-orientable, whence $\tilde{\mathbb{G}}$ is block-non-orientable □

Let us next consider the case of two adjacent components, that are not m-adjacent.

**Proposition 7.** *Suppose $\mathbb{G}$ is block-non-orientable and $C$ is a non-orientable component such that all orientable components adjacent to $C$ are not m-adjacent. Then there exists a m-ribbon $e$ in $C$, such that $\tilde{\mathbb{G}}$ obtained by slicing $e$ is block-non-orientable.*

*Proof.* Let $\mathbb{B}_C$ be the $\mathbb{G}$-block containing $C$. By Lemma 4, $C$ can sliced via some m-ribbon, $e$, into non-orientable components $C_1, \ldots, C_l$ and trivial components $T_1, \ldots, T_k$, where $T_1 = \{e_1\}, \ldots, T_k = \{e_k\}$.

First we show that each $\mathbb{G}$-block other than $\mathbb{B}_C$, is non-orientable in $\tilde{\mathbb{G}}$. Since slicing m-ribbon $e$ is defined as slicing with respect to the sectors $\gamma(e^L)$ and $e^R$, the path $Q^{\gamma(e^L), e^R}$ covers only $\mathbb{B}_C$. Lemma 6 guarantees that any $\mathbb{G}$-block other than $\mathbb{B}_C$ is a $\tilde{\mathbb{G}}$-block of the same orientability, i.e. it is non-orientable in $\tilde{\mathbb{G}}$.

Since $\mathbb{B}_C$ could be sliced into several $\tilde{\mathbb{G}}$-blocks, $\tilde{\mathbb{B}}_1, \ldots, \tilde{\mathbb{B}}_q$ and trivial components $T_1, \ldots, T_k$. It remains to show that each block $\tilde{\mathbb{B}}_i$ is non-orientable.

By construction, in each block $\tilde{\mathbb{B}}_i$, there exists a component $\tilde{C}_i$ adjacent to some trivial component $T_j$. We distinguish two cases:

1. $\tilde{C}_i$ is one of the $C_1, \ldots, C_l$, created by the slicing: as the $C_1, \ldots, C_l$ are all non-orientable, we are done.
2. $\tilde{C}_i$ corresponds to a $\mathbb{G}$-component in $\mathbb{B}_C$. Since $\tilde{C}_i$ is adjacent to a trivial component $T_j$ in $\tilde{\mathbb{G}}$, $\tilde{C}_i$ is adjacent to $C$ in $\mathbb{G}$. Thus the ribbon $e_j$ changes its crossing status and becomes the trivial component $T_j$ after slicing, which implies that $e_j$ intersects the interval $[\gamma(e^L), e^R]_\gamma$ in $\mathbb{G}$. Consequently, $e_j$ is a $m$-ribbon in $C$ and incident to the common vertex of $\tilde{C}_i$ and $C$. As a result $\tilde{C}_i$ and $C$ are $m$-adjacent. By assumption, all the orientable components adjacent to $C$ are not $m$-adjacent. Therefore $\tilde{C}_i$ and $\tilde{\mathbb{B}}_i$ are non-orientable.

$\square$

**Remark:** Suppose that $C$ is non-orientable and $C_1$ is orientable, adjacent but not $m$-adjacent to $C$. Let $r$ denote a $C$-ribbon incident to their common vertex and as such, $r$ is a $b$-ribbon. Successively applying Lemma 4 to $C$, i.e. successively slicing $m$-ribbons, renders $r$ as a trivial component at some point. This means that, along this process, $r$ changes its crossing status and Lemma 2 guarantees that then $r$ also changes its direction, becoming a $m$-ribbon.

Proposition 7 thus shows that after successive $C$-slicings, $C_1$ becomes in fact $m$-adjacent to a non-orientable component derived from successive slicings of $C$.

Fig. 13 shows how to transform an orientable component into non-orientable components by employing Propositions 6 and 7.

**Theorem 2.** *For each block-non-orientable unicellular fatgraph, $\mathbb{G}$ of genus $g$, we have $d(\mathbb{G}) = g$.*

*Proof.* We prove this by induction on $g$. Clearly, the assertion holds for $g = 1$. For $\mathbb{G}$ having genus $g$, we distinguish two cases:
Case 1: there exist one orientable component and one non-orientable component that are $m$-adjacent. Then Proposition 6 allows us to merge the two components by slicing and to obtain a block-non-orientable, unicellular fatgraph $\tilde{\mathbb{G}}$ of genus $g-1$.

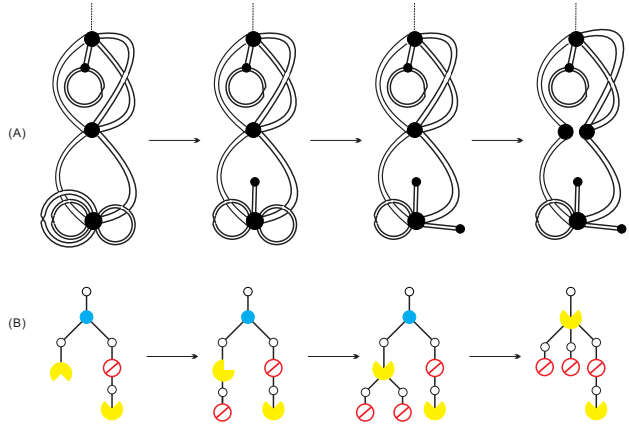

**Fig.** 13: The Pac-Man game: (A): via successive slicings, two components not $m$-adjacent first become $m$-adjacent and then merge into a non-orientable component. (B): the component trees: we represent a non-orientable component by a Pac-Man, an orientable component by a blue dot, and a trivial component by a no-entry sign.

Case 2: for any non-orientable component, all the orientable components are not $m$-adjacent. By Proposition 7, slicing an appropriate $m$-ribbon, $e$, of a non-orientable component we obtain a block-non-orientable, unicellular fatgraph $\tilde{\mathbb{G}}$ of genus $g-1$. □

## 9. A Hannenhalli-Pevzner formula for $r$-distance

First we study how half-flipping and gluing transforms orientable blocks into non-orientable blocks. Let $\mathbb{G}$ be a unicellular fatgraph having an orientable block $\mathbb{B}_1$. Let $C$ be an orientable component in $\mathbb{B}_1$ and $r$ be a $C$-ribbon. By $r^L, r^R$-half-flipping, we obtain a unicellular fatgraph $\tilde{\mathbb{G}}$.

**Proposition 8.** *The following assertions hold:*
*(1) each $\mathbb{G}$-block other than $\mathbb{B}_1$ is a $\tilde{\mathbb{G}}$-block of the same orientability and $\mathbb{B}_1$ becomes a non-orientable $\tilde{\mathbb{G}}$-block,*
*(2) the block tree $\tilde{B}$ of $\tilde{\mathbb{G}}$ is the same as the block tree $B$ of $\mathbb{G}$.*

*Proof.* Since $C$ is orientable, $r$ is a $b$-ribbon and sectors $r^L$ and $r^R$ are located at the same vertex and have the same orientation. Thus the $r^L, r^R$-half-flipping is well-defined. Note that the path $P^{r^L, r^R}$ contains only one component $C$ and thus the path $Q^{r^L, r^R}$ covers only $\mathbb{B}_1$. In view of Lemma 6, it remains to show that $\mathbb{B}_1$ is non-orientable in $\tilde{\mathbb{G}}$. Lemma 3 shows that $C$ is a $\tilde{\mathbb{G}}$-component. By construction, $r$ becomes a $m$-ribbon in $\tilde{\mathbb{G}}$, whence both $C$ and $\mathbb{B}_1$ are non-orientable in $\tilde{\mathbb{G}}$. By definition of block trees, the adjacency between any gap and any block or trivial component remains unchanged, proving assertion (2). $\square$

Let $\mathbb{B}_1$ and $\mathbb{B}_2$ be two orientable $\mathbb{G}$-blocks and let $Q$ denote the path connecting $\mathbb{B}_1$ and $\mathbb{B}_2$ in the block tree $B$.

**Proposition 9.** *There exists an $i, j$-gluing such that $\tilde{\mathbb{G}}$ satisfies*
*(1) each $\mathbb{G}$-block not covered by $Q$ is a $\tilde{\mathbb{G}}$-block having the same orientability and all blocks covered by $Q$, together with all trivial components contained in $Q$ merge into the non-orientable $\tilde{\mathbb{G}}$-block $\tilde{\mathbb{B}}_*$,*
*(2) the block tree $\tilde{B}$ of $\tilde{\mathbb{G}}$ is obtained from $B$ by contracting all blocks covered by $Q$ and trivial components contained in $Q$ into one black vertex.*

*Proof.* Let $C_k$ be an orientable component of $\mathbb{B}_k$ where $k = 1, 2$. Let $P$ denote the path joining $C_1$ and $C_2$ in the component tree $T$. Let $[s_k, t_k]_\gamma$ be the gap adjacent to $C_k$ on $P$ for $k = 1, 2$. We select $i \in I_{C_1} \setminus \{s_1, t_1\}$ and $j \in I_{C_2} \setminus \{s_2, t_2\}$. By construction, the path $P^{i,j}$ contains all the components on $P$. By construction, the path $Q^{i,j}$ joins $\mathbb{B}_1$ and $\mathbb{B}_2$, i.e. $Q^{i,j} = Q$. In order to prove (1), in view of Lemma 6, it suffices to show that $\tilde{\mathbb{B}}_*$ is non-orientable. According to Lemma 5, all components contained in $P$ merge into the $\tilde{\mathbb{G}}$-component $C_*$ in $\tilde{\mathbb{B}}_*$. Since at least one $b$-ribbon of $C_1$ intersects the interval $[i, j]_\gamma$ and changes to $m$-ribbon in $\tilde{\mathbb{G}}$, both $C_*$ and $\tilde{\mathbb{B}}_*$ are non-orientable. To prove (2), we verify, employing Lemma 1, that the adjacencies between gaps and blocks or trivial components in $\tilde{B}$ are induced from $B$: if a block not covered by $Q$ and a gap are adjacent in $B$, then they remain adjacent in $\tilde{B}$. If a gap is adjacent to a block covered by $Q$ then it is adjacent to $\tilde{\mathbb{B}}_*$ in $\tilde{B}$. Therefore the block tree $\tilde{B}$ of $\tilde{\mathbb{G}}$ is obtained from $B$ by contracting all blocks covered by $Q$ and trivial components on $Q$ into one block. $\square$

In difference to non-orientable blocks, that could be resolved individually, orientable blocks have to be considered as an ensemble. This is due to Proposition 9: the action of reversals affects the entire set of components or blocks.

**Definition 7. (E-block)** A block, $\mathbb{B}_i$, is *exposed* (*E-block*) if it is orientable and not covered by any path joining any two other orientable blocks in the block tree $B$. Let $h$ denote the number of E-blocks.

E-blocks are closely related to hurdles of signed permutations and $\pi$-maps [10, 13].

Let $\bar{B}$ denote the minimum $B$-subtree, containing all orientable blocks. By construction, each leaf of $\bar{B}$ is an orientable block. A block, $\mathbb{B}_i$, is an E-block if and only if $\mathbb{B}_i$ is a $\bar{B}$-leave and the gap adjacent to $\mathbb{B}_i$ has degree two in $\bar{B}$.

As a result, suppose an E-block, $\mathbb{B}_i$, is replaced by a non-orientable block, then, any other exposed block remains exposed and furthermore, at most one additional, orientable block becomes exposed.

**Definition 8. (S-block)** An E-block, $\mathbb{B}_i$, is called a *super-block* (S-block) if changing $\mathbb{B}_i$ into a non-orientable block does not affect the number of E-blocks in the block tree $B$.

**Proposition 10.** *Let $\tilde{\mathbb{G}}$ be the fatgraph obtained from $\mathbb{G}$ via an $i,j$-reversal. Then*

$$g(\mathbb{G}) + h(\mathbb{G}) - \big(g(\tilde{\mathbb{G}}) + h(\tilde{\mathbb{G}})\big) \leq 1. \tag{9.1}$$

*That is, a reversal reduces $(g+h)$ by at most one.*

*Proof.* We distinguish three cases

1. $i,j$-reversal is a slicing, i.e. $i$ and $j$ are located at the same vertex with different orientations. By construction, they are attached to components or gaps in the non-orientable block, $\mathbb{B}_0$. Lemma 6 shows that the $i,j$-reversal affects only $\mathbb{B}_0$ and that any E-block persists in $\tilde{\mathbb{G}}$. Thus the number of E-blocks of $\tilde{\mathbb{G}}$ satisfies $h(\tilde{\mathbb{G}}) \geq h(\mathbb{G})$. Since $g(\tilde{\mathbb{G}}) = g(\mathbb{G}) - 1$, eq. (9.1) holds.

2. $i,j$-reversal is half-flipping, i.e. $i$ and $j$ are located at the same vertex with the same orientation. By construction, they are attached to components or gaps in the block, $\mathbb{B}_1$. Lemma 6 shows that the $i,j$-reversal affects only $\mathbb{B}_1$. Proposition 8 shows that $\tilde{\mathbb{G}}$ has the same block tree as $\mathbb{G}$ and all $\mathbb{G}$-blocks except $\mathbb{B}_1$ have the same orientability. Furthermore, any E-block other than $\mathbb{B}_1$ persists in $\tilde{\mathbb{G}}$. Thus the number of E-blocks is reduced by at most one, i.e. $h(\tilde{\mathbb{G}}) \geq h(\mathbb{G}) - 1$ and in view of $g(\tilde{\mathbb{G}}) = g(\mathbb{G})$, eq. (9.1) follows.

3. $i,j$-reversal is a gluing, i.e. $i$ and $j$ are located at two distinct vertices. We have two cases

    - $i$ and $j$ are located at two distinct vertices of the same block, $\mathbb{B}_2$. Then Lemma 6 shows that all $\mathbb{G}$-blocks except $\mathbb{B}_2$ keep their

orientability in $\tilde{\mathbb{G}}$ and in particular, any E-block other than $\mathbb{B}_2$ persists to be an E-block in $\tilde{\mathbb{G}}$. Consequently, the number of E-blocks of $\tilde{\mathbb{G}}$ is reduced by at most one, i.e. $h(\tilde{\mathbb{G}}) \geq h(\mathbb{G}) - 1$, whence eq. (9.1).

- $i$ and $j$ are located at two distinct vertices in the two blocks, $\mathbb{B}_3$ and $\mathbb{B}_4$. Proposition 9 shows that the block tree, $\tilde{B}$, of $\tilde{\mathbb{G}}$ is obtained from $B$ by contracting all blocks covered by the path connecting $\mathbb{B}_3$ and $\mathbb{B}_4$. Thus any E-block other than $\mathbb{B}_3$ and $\mathbb{B}_4$ remains to be an E-block in $\tilde{\mathbb{G}}$. As a result, the number of E-blocks of $\tilde{\mathbb{G}}$ is reduced by at most two, i.e. $h(\tilde{\mathbb{G}}) \geq h(\mathbb{G}) - 2$. In view of $g(\tilde{\mathbb{G}}) = g(\mathbb{G}) + 1$, eq. (9.1) follows.

□

**Corollary 3.** *For any unicellular fatgraph $\mathbb{G}$, we have $d(\mathbb{G}) \geq g + h$.*

**Proposition 11.** *If the block tree $B$ has $2k$ E-blocks, then there exists a collection of $k$ B-paths joining all E-blocks such that each orientable block is covered by at least one path.*

**Theorem 3.** *Let $\mathbb{G}$ be a unicellular fatgraph of genus $g$ with $h$ E-blocks. Then*

$$d(\mathbb{G}) = \begin{cases} g+h+1 & \text{if } h \neq 1, h \text{ is odd and all } h \text{ E-blocks are super-blocks,} \\ g+h & \text{otherwise.} \end{cases}$$
(9.2)

*Proof.* Suppose first that $h \neq 1$, $h$ is odd and all $h$ E-blocks are S-blocks. Note that, the number of S-blocks cannot be reduced by half-flipping. Since gluing reduces the number of S-blocks by at most two, we need at least $(h+1)/2$ number of gluings to eliminate all S-blocks. Since each gluing increases the genus by one, we obtain the lower bound $g+(h+1)/2+(h+1)/2 = g+h+1$ on $r$-distance $d(\mathbb{G})$.

In view of Corollary 3, it is sufficient to identify a sequence of reversals that realizes the above lower bound. We distinguish the following cases

1. $h = 1$. Then $\mathbb{G}$ has only one orientable block and Proposition 8 shows that half-flipping one ribbon in this block transforms $\mathbb{G}$ into a block-non-orientable fatgraph $\mathbb{G}_1$ of genus $g$. Theorem 2 shows that $d(\mathbb{G}_1) = g$, whence $d(\mathbb{G}) = g + 1$.

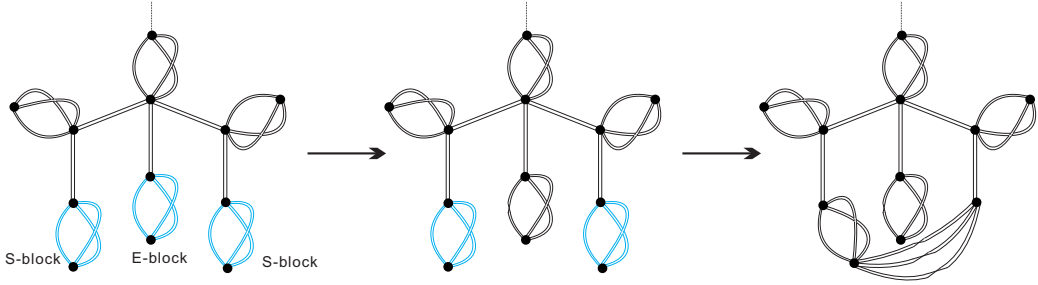

**Fig.** 14: Eliminating all the orientable blocks. The fatgraph $\mathbb{G}$ (Left) has genus 12 and contains 3 E-blocks, 2 of which are S-blocks. All orientable $\mathbb{G}$-blocks are transformed into non-orientable blocks in $\mathbb{G}'$ (Right), via one half-flipping and one slicing. Theorem 2 shows that $d(\mathbb{G}') = g(\mathbb{G}') = 13$. Therefore we have $d(\mathbb{G}) = g(\mathbb{G}) + h(\mathbb{G}) = 15$.

2. $h = 2k$. By Proposition 11, there exists a collection of $k$ paths joining all E-blocks and covering all orientable blocks. Proposition 9 implies, that one can always glue blocks covered by each path and thus obtain a block-non-orientable fatgraph, $\mathbb{G}_2$, via $k$ gluings. Therefore $d(\mathbb{G}_2) = g(\mathbb{G}_2) = g + k$ and $d(\mathbb{G}) = d(\mathbb{G}_2) + k = g + h$.
3. $h = 2k + 1$ and there exists at least one E-block, $\mathbb{B}_1$, which is not a S-block. Proposition 8 shows that half-flipping one $\mathbb{B}_1$-ribbon transforms $\mathbb{G}$ into a fatgraph $\mathbb{G}_3$ having genus $g$ and $2k$ E-blocks. As in the case of $h = 2k$, we have $d(\mathbb{G}_3) = g + 2k$ and thus $d(\mathbb{G}) = d(\mathbb{G}_3) + 1 = g + h$.
4. $h = 2k + 1$ and all $h$ E-blocks are S-blocks. Let $\mathbb{B}_1, \ldots, \mathbb{B}_h$ be all $h$ S-blocks. Set $B'$ to be the minimum subtree of $B$ containing all $\mathbb{B}_2, \ldots, \mathbb{B}_h$, and let $P'$ be the path connecting $\mathbb{B}_1$ and $B'$. By Proposition 9, we obtain the fatgraph $\mathbb{G}_4$ by gluing all blocks covered by $P'$. By construction, the number of E-blocks in $\mathbb{G}_4$ is reduced to $2k$. The case $h = 2k$ shows that $d(\mathbb{G}_4) = g + 1 + 2k$, whence $d(\mathbb{G}) = d(\mathbb{G}_4) + 1 = g + h + 1$.

□

Fig. 14 depicts a sequence of reversals eliminating E- and S-blocks.

**Acknowledgments**

The second author is a Thermo Fisher Scientific Fellow in Advanced Systems for Information Biology and acknowledges their support of this work.

[1] J. E. Andersen, R. C. Penner, C. M. Reidys, M. S. Waterman, Topological classification and enumeration of RNA structures by genus, J. Math. Biol. 67 (5) (2013) 1261–1278.

[2] D. Bessis, C. Itzykson, J. B. Zuber, Quantum field theory techniques in graphical enumeration, Adv. in Appl. Math. 1 (2) (1980) 109–157.

[3] M. A. Blasco, Telomeres and human disease: ageing, cancer and beyond, Nature Reviews Genetics 6 (8) (2005) 611–622.

[4] B. Bollobás, O. Riordan, A Polynomial Invariant of Graphs On Orientable Surfaces, Proc. Lond. Math. Soc. 83 (3) (2001) 513–531.

[5] G. Chapuy, A new combinatorial identity for unicellular maps, via a direct bijective approach, Adv. in Appl. Math. 47 (4) (2011) 874–893.

[6] G. Chapuy, V. Féray, E. Fusy, A simple model of trees for unicellular maps, J. Combin. Theory Ser. A 120 (8) (2013) 2064–2092.

[7] J. A. Ellis-Monaghan, I. Moffatt, Graphs on Surfaces: Dualities, Polynomials, and Knots, SpringerBriefs in Math., Springer-Verlag, New York, 2013.

[8] M. Esteller, Non-coding RNAs in human disease, Nature Reviews Genetics 12 (12) (2011) 861–874.

[9] A. Grothendieck, Esquisse d'un programme, in: L. Schneps, P. Lochak (Eds.), Geometric Galois Actions 1, Vol. 242 of London Math. Soc. Lecture Note Ser., Cambridge University Press, 1997, pp. 5–48.

[10] S. Hannenhalli, P. A. Pevzner, Transforming Cabbage into Turnip: Polynomial Algorithm for Sorting Signed Permutations by Reversals, J. ACM 46 (1) (1999) 1–27.

[11] J. Harer, D. Zagier, The Euler characteristic of the moduli space of curves, Invent. Math. 85 (1986) 457–485.

[12] F. W. D. Huang, C. M. Reidys, Topological language for RNA, Mathematical Biosciences 282 (2016) 109–120.

[13] F. W. D. Huang, C. M. Reidys, A topological framework for signed permutations, Discrete Math. 340 (9) (2017) 2161–2182.

[14] M. Kontsevich, Intersection theory on the moduli space of curves and the matrix airy function, Comm. Math. Phys. 147 (1) (1992) 1–23.

[15] S. K. Lando, A. K. Zvonkin, Graphs on Surfaces and Their Applications, Vol. 141 of Encyclopaedia Math. Sci., Springer-Verlag, Berlin, 2004, with an appendix by Don B. Zagier.

[16] B. Mohar, C. Thomassen, Graphs on Surfaces, Johns Hopkins Stud. Math. Sci., Johns Hopkins University Press, Baltimore, MD, 2001.

[17] H. Orland, A. Zee, RNA folding and large $n$ matrix theory, Nucl. Phys. B 620 (2002) 456–476.

[18] R. C. Penner, Perturbative series and the moduli space of Riemann surfaces, J. Differential Geom. 27 (1) (1988) 35–53.

[19] R. C. Penner, Moduli spaces and macromolecules, Bull. Amer. Math. Soc. 53 (2) (2016) 217–268.

[20] R. C. Penner, M. Knudsen, C. Wiuf, J. E. Andersen, Fatgraph models of proteins, Comm. Pure Appl. Math. 63 (10) (2010) 1249–1297.

[21] C. M. Reidys, F. W. D. Huang, J. E. Andersen, R. C. Penner, P. F. Stadler, M. E. Nebel, Topology and prediction of RNA pseudoknots. Bioinformatics 27 (8) (2011) 1076–1085.